\documentclass{amsart}

\usepackage{graphicx, pinlabel}
\usepackage{amssymb,amsmath,amsfonts}
\newtheorem{theorem}{Theorem}[section]
\newtheorem{proposition}[theorem]{Proposition}
\newtheorem{lemma}[theorem]{Lemma}
\newtheorem{conjecture}[theorem]{Conjecture}

\theoremstyle{definition}

\newtheorem{example}[theorem]{Example}

\newcommand{\Q}{\Bbb Q}
\newcommand{\R}{\Bbb R}
\newcommand{\G}{G(K)}
\newcommand{\C}{\Bbb C}
\newcommand{\Z}{\Bbb Z}
\newcommand{\D}{\Delta}

\newcommand{\tr}{{\mathrm{tr}\,}}
\newcommand{\la}{\langle}
\newcommand{\ra}{\rangle}
\newcommand{\p}{\partial}

\theoremstyle{remark}
\newtheorem{remark}[theorem]{Remark}
\numberwithin{equation}{section}

\begin{document}

\title[Twisted Alexander polynomials and character varieties]
{Twisted Alexander polynomials and character varieties of 2-bridge knot groups}

\author{Taehee Kim and Takayuki Morifuji}

\begin{abstract}
We study the twisted Alexander polynomial from the viewpoint of the
$SL(2,\C)$-character variety of nonabelian representations of a knot
group. It is known that if a knot is fibered, then the twisted
Alexander polynomials associated with nonabelian
$SL(2,\C)$-representations are all monic. In this paper, we show
that the converse holds for 2-bridge knots. Furthermore we show that
for a 2-bridge knot there exists a curve component in the
$SL(2,\C)$-character variety such that if the knot is not fibered
then there are only finitely many characters in the component for
which the associated twisted Alexander polynomials are monic.
We also show that for a 2-bridge knot of genus $g$, in the above
curve component for all but finitely many characters  the associated
twisted Alexander polynomials have degree $4g-2$.

\end{abstract}

\thanks{2000 {\it Mathematics Subject Classification}.
Primary 57M27, Secondary 57M05, 57M25.}

\thanks{{\it Key words and phrases.\/}
Twisted Alexander polynomial,
character variety, 2-bridge knot}

\address{Department of Mathematics, Konkuk University, Seoul 143-701,
Republic of Korea}
\email{tkim@konkuk.ac.kr}

\address{Department of Mathematics,
Tokyo University of Agriculture and Technology,
Tokyo 184-8588, Japan}
\email{morifuji@cc.tuat.ac.jp}

\maketitle

\section{Introduction}

In \cite{Lin01-1}, Lin introduced the twisted Alexander polynomial
for knots in the $3$-sphere by using regular Seifert surfaces, while
Wada \cite{Wada94-1} defined it for finitely presentable groups,
which include link groups. It is a generalization of the classical
Alexander polynomial and has many applications to knot theory. A
notable application is to distinguish two mutant knots with the
trivial Alexander polynomial \cite{Wada94-1}. The twisted Alexander
polynomial can be regarded as the Reidemeister torsion and this fact
leads to its symmetry \cite{Kitano96-1}. It is also useful to the
problems on inversion and concordance of knots \cite{KL99-1}.
Furthermore classical results about the Alexander polynomial of
fibered knots are extended to the twisted case
\cite{Cha03-1,FK05-1,FV08-1,GKM05-1}. In particular, in
\cite{FV08-1} Friedl and Vidussi showed that the twisted Alexander
polynomials corresponding to all finite representations detect
fibered $3$-manifolds. For literatures on the twisted Alexander
polynomial and other related topics, refer to the survey paper by
Friedl and Vidussi \cite{FV08-2}.

The purpose of this paper is to consider another approach to the
fibering problem for knots. That is, we study the problem from the
viewpoint of the $SL(2,\C)$-character variety of a knot group. The
approach of using the character variety to 3-manifold problems was
introduced by Culler and Shalen \cite{CS83-1}, and since then, many
deep results on 3-manifolds have been obtained using the character
variety. Roughly speaking, many topological properties of a
$3$-manifold are encoded in the character variety of the 3-manifold
group, and therefore it is interesting to study basic properties of
the character variety.

In general, for a given knot each coefficient of the twisted
Alexander polynomial of the knot defines a complex valued function
on the $SL(2,\C)$-representation variety of the knot group. If the
coefficient of the highest degree term is the constant function 1,
then we call such a representation (and its character)
\textit{monic}.\/ In this point of view, it is known that every
nonabelian $SL(2,\C)$-representation of a fibered knot is monic
\cite{GKM05-1}. One can ask if the converse holds, but the answer is
still unknown (cf. \cite{FV08-1}). In this paper, we give a partial
answer that the converse holds for 2-bridge knots: for a 2-bridge
knot, if every nonabelian $SL(2,\C)$-representation of the knot is
monic, then the knot is fibered (Theorem~\ref{thm:fibered-monic}).
In fact, we show a much stronger theorem that for a nonfibered
2-bridge knot, there exists a curve component in the character
variety of nonabelian representations of the knot in which there are
only finitely many monic characters
(Theorem~\ref{thm:fibered-finite-1}). This extends the result for
twist knots of the second author in \cite{Morifuji08-1}. We also
give a sufficient condition based on the Alexander polynomial for a
knot to have only finitely many monic characters
(Theorem~\ref{thm:fibered-finite-2}). Moreover we give an upper
bound of the number of monic characters for a family of certain
2-bridge knots which contains twist knots (Theorem~\ref{thm:5.1}).
Although fiberedness of a 2-bridge knot is detected by the Alexander
polynomial since 2-bridge knots are alternating, these results give
new fibering criteria for 2-bridge knots and lead us to a conjecture
that Theorem~\ref{thm:fibered-finite-1} can be extended for a
general nonfibered knot (see Conjecture~\ref{conj:6.4} and
Remark~\ref{rmk:6.5}).

In addition to studying the fibering problem, we also investigate
detecting the genus of a knot using the twisted Alexander
polynomial. Recall that a 2-bridge knot is alternating, hence its
classical Alexander polynomial is of degree $2g$ where $g$ denotes
the genus of the knot. It is also known that if a knot of genus $g$
is fibered, then the twisted Alexander polynomial associated with
any nonabelian $SL(2,\C)$-representation is of degree $4g-2$ (see
\cite{KM05-1}). In this direction, we show that for a (possibly
nonfibered) 2-bridge knot of genus $g$, there exists a curve
component in the character variety of nonabelian representations of
the knot such that for all but finitely many characters in the
component the associated twisted Alexander polynomials have degree
$4g-2$ (Theorem~\ref{thm:genus1}).

This paper is organized as follows. In the next section, we review
some basic materials for the character variety, in particular, the
representation polynomial of a 2-bridge knot group due to Riley
\cite{Riley84-1}. In Section 3, we quickly review the definition of
the twisted Alexander polynomials associated with
$SL(2,\C)$-representations. In Section 4, we show finiteness of
monic characters in a curve component of the character variety of
nonabelian representations of a nonfibered 2-bridge knot. 
We also show the genus of a
2-bridge knot is detected by the twisted Alexander polynomial and
the character variety of the knot group.
Section 5 is devoted to the calculation of an upper bound of the number of
monic characters for a family of 2-bridge knots. In the last
section, we give some remarks and state conjectures on the fibering
problem from the viewpoint of $SL(2,\C)$-representations of a knot
group.

\section{Character varieties}
\label{sec:chracter varieties} In this section we review the theory
of the character variety of a knot group that will be needed for our
purpose. See \cite{CS83-1,Le93-1,Riley84-1} for details.

Let $G$ be a finitely generated group. Define
$R(G)=\mathrm{Hom}(G,SL(2,\C))$ to be the set of representations of
$G$ into $SL(2,\C)$. It is known that $R(G)$ is a complex affine
algebraic set. We call it the $SL(2,\C)$-{\it representation
variety}\/ of $G$, though it might be a union of a finite number of
irreducible algebraic varieties in the sense of algebraic geometry.
The isomorphism class of this variety does not depend on the choice
of generators.

The group $SL(2,\C)$ acts on $R(G)$ by conjugation. Let $\hat{R}(G)$
denote the set of orbits. Two representations $\rho,\rho'\in R(G)$
are called \textit{conjugate}\/ if they lie in the same orbit. The
algebro-geometric quotient of $R(G)$ under this conjugate action is
called the $SL(2,\C)$-\textit{character variety}\/ of $G$, which is
denoted by $X(G)$. The \textit{character}\/ of a representation
$\rho$ is a map $\chi_\rho:G\to\C$ defined by
$\chi_\rho(\gamma)=\tr(\rho(\gamma))$ for $\gamma\in G$. There is a
bijection between $X(G)$ and the set of characters of $G$. Namely,
we have a canonical identification $X(G)=\{\chi_\rho\,|\,\rho\in
R(G)\}$.

A representation $\rho:G\to SL(2,\C)$ is said to be
\textit{abelian}\/ if $\rho(G)$ is an abelian subgroup of
$SL(2,\C)$. A representation $\rho$ is called \textit{reducible}\/
if there exists a proper invariant subspace in $\C^2$. This is
equivalent to saying that $\rho$ can be conjugated to a
representation by upper triangular matrices. It is easy to see that
every abelian representation is reducible, but the converse does not
hold. Namely there is a reducible nonabelian representation in
general. When $\rho$ is not reducible, it is called
\textit{irreducible}.\/ If $\rho,\rho'\in R(G)$ have the same
character $\chi_\rho=\chi_{\rho'}$ and $\rho$ is irreducible, then
$\rho$ is conjugate to $\rho'$ (see \cite[Proposition
1.5.2]{CS83-1}).

Let $R^{\mathrm{irr}}(G)$ denote the subset of irreducible
representations of $G$ and $X^{\mathrm{irr}}(G)$ denote its image
under the map $t:R(G)\to X(G)$ given by $t(\rho)=\chi_\rho$.
Similary, we write $X^{\mathrm{nab}}(G)$ for the image of
$R^{\mathrm{nab}}(G)$, the set of nonabelian
$SL(2,\C)$-representations, by $t$.

If $G$ is a knot group $\G$,
namely the fundamental group of the exterior $E(K)$ of
a knot $K$ in the 3-sphere $S^3$,
we denote $R(\G)$ and $X(\G)$ by
$R(K)$ and $X(K)$ for simplicity.
We also use the notations
$X^{\mathrm{irr}}(K),\,X^{\mathrm{nab}}(K)$ and so on.

Now
we quickly review how to describe the nonabelian part of
the $SL(2,\C)$-character variety of a $2$-bridge knot group
(see \cite{Le93-1,Riley84-1} for details).
Let
$K=K(\alpha,\beta)$ be a 2-bridge knot corresponding to
a pair of relatively prime odd integers $(\alpha,\beta)$
with $-\alpha<\beta<\alpha$.
The knots
$K(\alpha,\beta)$ and $K(\alpha',\beta')$ have the same type
if and only if
$\alpha=\alpha'$ and $\beta\equiv \beta'$ or
$\beta\beta'\equiv 1\mod \alpha$.
The knot group of $K$ has a presentation
$$
\G
=
\la
a,b\,|\,wa=bw
\ra,\quad
w=a^{\epsilon_1}b^{\epsilon_2}\cdots
a^{\epsilon_{\alpha-2}}b^{\epsilon_{\alpha-1}}
$$
where
$\epsilon_i=(-1)^{[\frac{\beta}{\alpha}i]}$
and
$[u]$ denotes the greatest integer less than or
equal to $u\in\R$.
It follows that
$\epsilon_i=\epsilon_{\alpha-i}$ holds for any $i$.

The above presentation for $G(K)$ is not unique,
but
the existence of at least one such presentation
follows from
Wirtinger's algorithm applied to Schubert's canonical
$2$-bridge diagram of $K(\alpha,\beta)$.
The generators $a$ and $b$ come from the two bridges
(overpasses) and represent the meridian up to
conjugation.

We consider the matrices
$C=
\begin{pmatrix}
s & 1 \\
0 & s^{-1}
\end{pmatrix}
$
and
$D=
\begin{pmatrix}
s & 0 \\
2-y & s^{-1}
\end{pmatrix}
$, where $s\not=0,\, y\in\C$. Since $a$ and $b$ are conjugate in the
group $\G$, for a nonabelian representation $\rho:\G\to SL(2,\C)$,
$\rho(a)$ and $\rho(b)$ have the same trace. Moreover, taking
conjugation if necessary, we can assume $\rho(a)=C$ and $\rho(b)=D$
without loss of generality. Here the entry $2-y$ is chosen so that
$CD^{-1}$ has trace $y$ (see \cite{MPV09-1}). Under the setting
above, we have the following.

\begin{proposition}
\cite[Theorem 1]{Riley84-1} \label{pro:2.1} The assignment
$\rho(a)=C,\,\rho(b)=D$ defines a nonabelian representation of $\G$
if and only if the pair $(s,y)$ satisfies the equation
$$
w^{11}+(s^{-1}-s)w^{12}=0,
$$
where
$W=\rho(w)=(w^{ij})$.
Conversely,
every nonabelian representation of $\G$ is conjugate to
a representation satisfying the above equation.
\end{proposition}

We now define the polynomial $\phi(s,y)$ via
$$
\phi(s,y)
=
w^{11}+(s^{-1}-s)w^{12}\in {\Z}[s^{\pm1},y]
$$
and
call it the \textit{Riley polynomial}\/ of a 2-bridge knot $K$.
That is,
the Riley polynomial gives a defining equation of the nonabelian part
$\hat{R}^{\mathrm{nab}}(K)$ of
conjugacy classes of $SL(2,\C)$-representations.
We also remark here that
$\phi(s,y)$ might be a reducible polynomial
even over ${\Z}[s^{\pm1},y]$.

\begin{example}\label{ex:2.3}
Let $K$ be the knot $K(15,11)=7_4$ with
the Alexander polynomial $\D_K(t)=4-7t+4t^2$.
The knot group has a presentation
$$
G(K)
=
\la
a,b\,|\, wa=bw
\ra,
\quad
w=\left((ba^{-1})^2(b^{-1}a)^2\right)^2.
$$
A direct calculation shows that
the Riley polynomial $\phi(s,y)$ has the
factorization
$$
\phi(s,y)=\phi_1(s,y)\phi_2(s,y)\in{\Z}[s^{\pm1},y],
$$
where
\begin{align*}
\phi_1(s,y)
&=
1+(s^{-2}+s^2)y^2-y^3,
\\
\phi_2(s,y)
&=
1
-2(s^{-2}+s^2)y
+(3s^{-2}+2+3s^2)y^2
-(s^{-2}+3+s^2)y^3
+y^4.
\end{align*}
Putting $y=2$,
$\phi_1(s,2)=4s^{-2}-7+4s^2=s^{-2}\D_K(s^2)=0$ gives reducible
nonabelian representations
(see Proposition \ref{pro:2.4} below).
On the other hand,
$\phi_2(s,y)=0$ contains no reducible representation.
\end{example}

For $\gamma\in G(K)$, let $t_\gamma$ be a function $t_\gamma:R(K)\to
\C$ defined by $t_\gamma(\rho)=\tr(\rho(\gamma))$. Then as a
coordinate of $X^{\mathrm{nab}}(K)$ for a 2-bridge knot $K$, we can
take $(t_a,t_{ab^{-1}})$. Namely, $X^{\mathrm{nab}}(K)$ can be
identified with the image of $R^{\mathrm{nab}}(K)$ under the map
$(t_a,t_{ab^{-1}}):R(K)\to \C^2$ (see
\cite[Proposition 1.4.1]{CS83-1} and \cite[Section 2]{MPV09-1}). Of
course, this map factors through $\hat{R}(K)$. If $\phi$ is viewed
as a polynomial in $x=s+s^{-1}$ and $y$, then $X^{\mathrm{nab}}(K)$
is given in $\C^2$ by $\phi(x,y)=0$
(see \cite[Proposition 2.2]{MPV09-1}).

\begin{remark}\label{rmk:2.3}
It is known that
for a 2-bridge knot $K$
the Zariski closure $\overline{X^{\mathrm{irr}}(K)}$ is
exactly the nonabelian part of
the character variety $X^{\mathrm{nab}}(K)$
(see \cite{Le93-1}).
More precisely,
except finitely many cases,
a nonabelian representation is irreducible.
\end{remark}

We end this section with the following useful
proposition.

\begin{proposition}
[Burde \cite{Burde67-1}, de Rham \cite{deRham67-1}]
\label{pro:2.4}
Let
$\eta_0:G(K)\to SL(2,\C)$ be an abelian representation of a knot $K$
given by
$\eta_0(\mu)
=
\begin{pmatrix}
\lambda & 0 \\
0 & \lambda^{-1}
\end{pmatrix}
$,
where
$\mu$ is the meridian of $K$ and $\lambda\not=0\in\C$.
Then
there is a reducible nonabelian representation
$\rho:G(K)\to SL(2,\C)$ so that
$\chi_\rho=\chi_{\eta_0}$
if and only if
$\D_K(\lambda^2)=0$.
\end{proposition}

\section{Twisted Alexander polynomials}
\label{sec:tap}

For a knot group $\G=\pi_1(E(K))$, we choose and fix a Wirtinger
presentation
$$
\G=
\langle \gamma_1,\ldots,\gamma_k\,|\,r_1,\ldots,r_{k-1}\rangle.
$$
Then
the abelianization homomorphism
$$
\alpha:\G\to H_1(E(K),{\Bbb Z})
\cong {\Z}
=\la t
\ra
$$
is given by $\alpha(\gamma_1)=\cdots=\alpha(\gamma_k)=t$. Here we
specify a generator $t$ of $H_1(E(K),\Z)$ and denote the sum in $\Z$
multiplicatively. 
In this paper we consider just a linear representation
$\rho:\G\to GL(2,\C)$.

These maps $\rho$ and $\alpha$ naturally induce two ring
homomorphisms $\tilde{\rho}: {\Z}[G(K)] \rightarrow M(2,{\C})$ and
$\tilde{\alpha}:{\Z}[G(K)]\rightarrow {\Z}[t^{\pm1}]$, where
${\Z}[G(K)]$ is the group ring of $G(K)$ and $M(2,{\C})$ is the
matrix algebra of degree $2$ over ${\C}$. Then
$\tilde{\rho}\otimes\tilde{\alpha}$ defines a ring homomorphism
${\Z}[G(K)]\to M\left(2,{\C}[t^{\pm1}]\right)$. Let $F_k$ denote the
free group on generators $\gamma_1,\ldots,\gamma_k$ and
$$
\Phi:{\Z}[F_k]\to M\left(2,{\C}[t^{\pm1}]\right)
$$
the composition of the surjection
$p:{\Z}[F_k]\to{\Z}[G(K)]$
induced by the presentation of $G(K)$
and the map
$\tilde{\rho}\otimes\tilde{\alpha}:{\Z}[G(K)]\to M(2,{\C}[t^{\pm1}])$.

Let us consider the $(k-1)\times k$ matrix $M$
whose $(i,j)$-entry is the $2\times 2$ matrix
$$
\Phi\left(\frac{\partial r_i}{\partial \gamma_j}\right)
\in M\left(2,{\C}[t^{\pm1}]\right),
$$
where
$\frac{\partial}{\partial \gamma}$
denotes the Fox differential.
This matrix $M$ is called the
\textit{Alexander matrix}\/ of
$\G$
associated with the representation $\rho$.

For
$1\leq j\leq k$,
let us denote by $M_j$
the $(k-1)\times(k-1)$ matrix obtained from $M$
by removing the $j$th column.
We regard $M_j$ as
a $2(k-1)\times 2(k-1)$ matrix with coefficients in
${\C}[t^{\pm1}]$.

Then Wada's \textit{twisted Alexander polynomial}\/ of a knot $K$
associated with a representation $\rho:\G\to GL(2,{\C})$ is defined
to be the rational function
$$
\D_{K,\rho}(t)
=\frac{\det M_j}{\det\Phi(1-\gamma_j)}
$$
and well-defined up to multiplication by 
$\varepsilon t^{2i}~(\varepsilon\in{\C}^*, i\in{\Z})$. 
Moreover for the case of a
special linear representation $\rho:\G\to SL(2, \C)$,
$\D_{K,\rho}(t)$ is well-defined up to multiplication by
$t^{2i}~(i\in{\Z})$. In this paper mostly we use a representation
$\rho\colon \G \to SL(2, \C)$.

\begin{remark}\label{rmk:3.1}
If $\rho$ is conjugate to $\rho'$ in $GL(2,\C)$, 
then $\D_{K,\rho}(t)=\D_{K,\rho'}(t)$ holds (see \cite[Section 3]{Wada94-1}).
Moreover $\D_{K,\rho}(t)$ has the following properties. 
In this paper, a (Laurent) polynomial
$f(t)=c_mt^m+c_{m-1}t^{m-1}+\cdots+c_{n+1}t^{n+1}+c_nt^n\in\C[t^{\pm
1}]$ is called \textit{monic}\/ if the coefficient $c_m$ is just 1
(not $\pm1$). Let $\rho\colon G(K)\to SL(2,\C)$ be a 
nonabelian representation.

\begin{enumerate}
\item
The twisted Alexander polynomial $\D_{K,\rho}(t)$  is always a
polynomial for any knot $K$ (see \cite[Theorem 3.1]{KM05-1}), and it
is reciprocal, i.e., $\D_{K,\rho}(t)=t^i\D_{K,\rho}(t^{-1})$ for
some $i\in \Z$ (see \cite[Corollary 3.5]{HSW10-1}).

\item
If $K$ is a \textit{fibered knot}\/ of genus $g$, namely the
exterior $E(K)$ has the structure of a surface bundle over the
circle, then $\D_{K,\rho}(t)$ becomes a monic polynomial of degree
$4g-2$ (see \cite[Theorem 3.1]{GKM05-1}  and \cite[Theorem
3.2]{KM05-1}).

\item 
If $K$ is a knot of genus $g$, then 
$\deg (\D_{K,\rho}(t))\le 4g-2$ 
(see \cite[Theorem 1.1]{FK05-1}).

\item 
If $\rho$ is a reducible nonabelian representation, then 
up to conjugation,
$\rho(\gamma_i)=
\begin{pmatrix}
\lambda & \nu_i \\
0&\lambda^{-1}
\end{pmatrix}$ for each $i$ where $\lambda\not=0,\,\nu_i\in\C$ and 
$$
\D_{K,\rho}(t) =\frac{\D_K(\lambda t)\D_K(\lambda^{-1}t)}
{(t-\lambda)(t-\lambda^{-1})}.
$$
(See the proof of \cite[Theorem 3.1]{KM05-1}.) In particular,
$\D_{K,\rho}(t)$ is a polynomial of degree $2\deg(\D_K(t))-2$.

\end{enumerate}
\end{remark}

\begin{example}\label{ex:3.2}
Let $\eta:G(K)\to SL(2,\C)$ be an abelian representation
defined by the correspondence
$\gamma_i\mapsto
\begin{pmatrix}
\lambda & \nu \\
0&\lambda^{-1}
\end{pmatrix}$ for any $i$,
where $\lambda\not=0,\,\nu\in\C$. 
Then
$$
\D_{K,\eta}(t) =\frac{\D_K(\lambda t)\D_K(\lambda^{-1}t)}
{(t-\lambda)(t-\lambda^{-1})}
$$
but it is not a polynomial in general.
\end{example}

Let $\rho:\G\to SL(2,\C)$ be a nonabelian representation of a knot
$K$. Then each coefficient of $\D_{K,\rho}(t)$ defines a complex
valued function on the representation variety $R^{\mathrm{nab}}(K)$,
the orbit space $\hat{R}^{\mathrm{nab}}(K)$ and also on the
character variety $X^{\mathrm{nab}}(K)$. In fact, if $\rho$ is a
reducible nonabelian representation, there is an abelian
representation $\eta_0$ so that $\chi_{\eta_0}=\chi_\rho$ (see
Proposition \ref{pro:2.4}). Of course, $\rho$ and $\eta_0$ are not
conjugate, but $\D_{K,\rho}(t)=\D_{K,\eta_0}(t)$ holds (see 
Remark~\ref{rmk:3.1} (iv) and Example \ref{ex:3.2}). Therefore each
coefficient of $\D_{K,\rho}(t)$ can be considered as a function on
$X^{\mathrm{nab}}(K)$. Also we define {\it the
twisted Alexander polynomial associated with $\chi\in
X^{\mathrm{nab}}(K)$} to be $\D_{K,\rho}(t)$ where $\chi=\chi_\rho$,
and denote it by $\D_{K,\chi}(t)$. In particular, 
for a 2-bridge knot $K$, since $X^\mathrm{nab}(K)$
is given in $\C^2$ by the equation $\phi(x,y)=0$ where $x=s+s^{-1}$, 
each coefficient of
$\D_{K,\rho}(t)$ can be considered as a function of $s$ and $y$ or a
function of $x$ and $y$.

It should be noted that the notion of monic
polynomial makes sense for $\D_{K,\rho}(t)$ by 
the indeterminacy of $\D_{K,\rho}(t)$ and
Remark \ref{rmk:3.1} (i). We say a 
nonabelian representation $\rho:G(K)\to SL(2,\C)$ is
\textit{monic} if the twisted Alexander polynomial $\D_{K,\rho}(t)$
associated with $\rho$ is a monic polynomial. 
Similary we call 
$\chi \in X^{\mathrm{nab}}(K)$ \textit{monic}\/ if $\D_{K,\chi}(t)$
is a monic polynomial. 
For a fibered knot $K$ in $S^3$, all the nonabelian
$SL(2,\C)$-representations are monic \cite{GKM05-1}. In other words,
the coefficient of the highest degree term of $\D_{K,\rho}(t)$ is a
constant function 1 on the whole $X^{\mathrm{nab}}(K)$.

\section{Finiteness theorems for $2$-bridge knots}

Below we show that for a 
2-bridge knot $K$ in $S^3$,
$X^{\mathrm{nab}}(K)$ detects if $K$ is fibered. That is, the
converse of Remark~\ref{rmk:3.1}~(ii) holds for 2-bridge knots.

\begin{theorem}\label{thm:fibered-monic}
A $2$-bridge knot $K$ is fibered if and only if $\Delta_{K,\rho}(t)$
is monic for any nonabelian representation $\rho:G(K)\to SL(2,\C)$.
\end{theorem}

\begin{proof}
The `if' part was proven in
\cite{GKM05-1}. Now assume that $K$ is a nonfibered 2-bridge knot.
We choose and fix a Wirtinger presentation $G(K)=\la
a,b\,|\,wa=bw\ra$ as we did in Section~\ref{sec:chracter varieties}.
Let $\rho:\G\to SL(2,\C)$ be an arbitrary nonabelian representation.
We will show that the coefficient of the highest degree term of
$\D_{K,\rho}(t)$ is not identically one on $X^{\mathrm{nab}}(K)$.

By taking conjugations if necessary, we may assume that
$$
\rho(a) =\begin{pmatrix}
s & 1 \\
0 & s^{-1}
\end{pmatrix},\qquad
\rho(b) =\begin{pmatrix}
s & 0 \\
2-y & s^{-1}
\end{pmatrix}.
$$
Suppose that the top coefficient
of $\D_{K,\rho}(t)$
is identically one on $X^{\mathrm{nab}}(K)$. Namely, we assume that
$$
\D_{K,\rho}(t) = t^m+f_1(s,y)t^{m-1}+f_2(s,y)t^{m-2}+\cdots
$$
holds. Since the Alexander polynomial of a 2-bridge knot is not
trivial, by Proposition~\ref{pro:2.4} there exists a reducible
nonabelian representation of $K$. Since $\rho$ becomes reducible
exactly when $y=2$ (see \cite[Proposition 1.5.5]{CS83-1}), it means
that there exists $s_0\in \C^*$ such that $\phi(s_0 + s_0^{-1},2) =
0$ viewing $\phi$ as a polynomial of $x=s+s^{-1}$ and $y$, and
therefore $\rho$ with $s=s_0$ and $y=2$ becomes a reducible
nonabelian representation, say $\rho'$. For $\rho'$, 
by Remark~\ref{rmk:3.1} (iv) we
obtain that $\D_{K,\rho'}(t) = \det M/(t-s)(t-s^{-1})$ where
\begin{equation}\label{eqn:det}
\det M = \D_K(st)\D_K(s^{-1}t).
\end{equation}
Here we remark that the top coefficient of $\D_{K,\rho'}(t)$ remains
the constant 1 by our assumption. Since $K$ is a nonfibered
$2$-bridge knot (in particular a nonfibered alternating knot), the
Alexander polynomial $\D_K(t)$ is nonmonic (see \cite[Theorem 1.2]{Murasugi63-1}).
Hence the right hand side of (\ref{eqn:det})
is also nonmonic, but this is a contradiction because
$\D_{K,\rho'}(t)$ is monic if and only if the numerator of
$\D_{K,\rho'}(t)$, which is $\det M$, is monic.
\end{proof}

The above theorem can be strengthened further. Namely, for a
nonfibered 2-bridge knot in $S^3$, we can find a special curve
component of the character variety which contains only finitely many
monic characters:

\begin{theorem}\label{thm:fibered-finite-1}
For a nonfibered $2$-bridge knot $K$, there exists an irreducible
curve component in $X^{\mathrm{nab}}(K)$ which contains only a
finite number of monic characters.
\end{theorem}

\begin{proof}
As we saw in the proof of Theorem~\ref{thm:fibered-monic}, the Riley
polynomial $\phi(x,y)$ of $K$ has an irreducible factor
$\phi_1(x,y)\in \C[x,y]$ such that  $\phi_1(s_0+s_0^{-1}, 2)=0$ for
some $s_0\in\C^*$. 
Let $X_1$ be the curve given by $\phi_1(x,y)=0.$
In particular, $X_1$ contains a character of a reducible nonabelian
representation of $K$. By Remark~\ref{rmk:3.1} (iii) we may assume
that for $\chi=\chi(x,y)\in X_1$,
$$
\D_{K,\chi}(t) = \psi_{4g-2}(x,y)t^{4g-2} +
\psi_{4g-3}(x,y)t^{4g-3}+\cdots +\psi_1(x,y)t +\psi_0(x,y).
$$
Then $\D_{K,\chi}(t)$ is monic if $\psi_{4g-2}(x,y)=1$ or there exists 
$j\ge 2g-1$ such that $\psi_i(x,y)=0$ for $j< i\le 4g-2$ and
$\psi_j(x,y)=1$. Therefore the set of monic characters in the
irreducible curve component $X_1$ is contained in $\{(x,y)\in
\C^2\,\mid\, \phi_1(x,y)=0 \mbox{ and } \psi_{4g-2}(x,y)=1\} \cup
\{(x,y)\in \C^2\,\mid\, \phi_1(x,y)=0 \mbox{ and }
\psi_{4g-2}(x,y)=0\}$. By our choice of $\phi_1$, $\psi_{4g-2}(x,y)$
is not identically one on $X_1$ as 
was shown in the proof of Theorem~\ref{thm:fibered-monic}. Also
since $K$ is alternating, $\deg(\D_K(t))= 2g$ 
(see \cite{Crowell59-1,Murasugi58-1}). Therefore by our
choice of $\phi_1$ and Remark~\ref{rmk:3.1} (iv), $\psi_{4g-2}(x,y)$
is not identically 0 on $X_1$. 
Now we can conclude that $\phi_1(x,y)$ and $\psi_{4g-2}(x,y)-1$ have
no common divisors, and the affine plane curves defined by them have
at most finitely many intersection points (B\'{e}zout's Theorem).
Similarly, $\phi_1(x,y)=0$ and $\psi_{4g-2}(x,y)=0$ have at most
finitely many intersection points. Therefore the number of monic
characters on $X_1$ is finite.
\end{proof}

Note that in the above proof we
do not need non-fiberedness of a 2-bridge knot to show that
$\psi_{4g-2}(x,y)$ is not identically 0. 
Therefore $\phi_1(x,y)=0$ and $\psi_{4g-2}(x,y)=0$ have at most
finitely many intersection points for {\it any} 2-bridge knot, and
we obtain the following theorem:

\begin{theorem}\label{thm:genus1}
For a $2$-bridge knot $K$ of genus $g$, there exists an irreducible
curve component $X_1$ in $X^{\mathrm{nab}}(K)$ such that
$\deg(\D_{K,\chi}(t)) = 4g-2$ for all but finitely many $\chi\in
X_1$.
\end{theorem}

\begin{remark}\label{rmk:4.3.1}
\begin{enumerate}
\item
It is known that for any positive integer $n$ there is a hyperbolic $2$-bridge knot $K$
such that $X^{\mathrm{irr}}(K)$ has at least $n$ irreducible components
(see \cite[Corollary 7.3]{ORS08-1}).

\item
For a knot $K$ of genus $g$, in \cite{DFJ10-1} Dunfield, Friedl and
Jackson show that the set of monic characters in $X(K)$ is Zariski
closed and $\{\chi\in X(K)\,\mid \, \deg(\D_{K,\chi}(t)) = 4g-2\}$
is Zariski open.
\end{enumerate}
\end{remark}

Assuming a suitable condition for the 2-bridge knot $K(\alpha,\beta)$,
we obtain the following finiteness result.

\begin{theorem}\label{thm:fibered-finite-2}
Let $K=K(\alpha,\beta)$ be a nonfibered $2$-bridge knot and $c\in\Z$
the leading coefficient of $\D_K(t)$. Let $p$ be an odd prime
divisor of $\alpha$. 
Suppose $c\not\equiv 0 \mbox{ and
} c^2\not\equiv \pm 1\,\, \mathrm{mod}~p$. Then the number of monic
characters in $X^\mathrm{nab}(K)$ is finite.
\end{theorem}

\begin{proof}
By Lemmas \ref{lem:metabelian} and \ref{lem:nonmonic} below, we see
that the coefficient of the highest degree term of the twisted
Alexander polynomials is not identically one on each irreducible
component of $X^{\mathrm{nab}}(K)$. Therefore, by the arguments in
the proof of Theorem~\ref{thm:fibered-finite-1} the number of monic
characters is finite for each irreducible component of
$X^\mathrm{nab}(K)$.
\end{proof}

A representation $\rho:G(K)\to SL(2,\C)$ is called
\textit{metabelian}\/ if the image of the commutator subgroup
$\left[G(K),G(K)\right]$ is an abelian subgroup in $SL(2,\C)$. It
should be noted that reducible representations are metabelian by
\cite[Lemma 1.2.1]{CS83-1} (see also \cite[Section 2]{NY08-1}).

\begin{lemma}\label{lem:metabelian}
Any irreducible component of
$X^\mathrm{nab}(K)$ of a $2$-bridge knot $K$ contains a character of
an irreducible metabelian representation.
\end{lemma}

\begin{proof}
As a polynomial in $\Z[s^{\pm1}][y]$, the coefficient of the highest
degree term of the Riley polynomial $\phi(s,y)$ is $\pm1$ (see
\cite[p.197]{Riley84-1}). Therefore $\phi(s,y)$ has the irreducible
factorization
$$
\phi(s,y)=\pm\phi_1(s,y)\phi_2(s,y)\cdots\phi_m(s,y)
$$
over $\Z[s^{\pm1}][y]$ such that each $\phi_k(s,y)$ is a monic
polynomial in $\Z[s^{\pm1}][y]$. By Lemma~3 and Proposition~1 in
\cite{Riley84-1}, it follows that $\phi_k(s,y)=\phi_k(s^{-1},y)$ for
any $k$. We then put $x=s+s^{-1}$ and regard each $\phi_k$ as an
element of $\Z[x,y]$. It might be a reducible polynomial over
$\C[x,y]$, but we see that $\phi_k(x,y)$ can be written as the
product of monic polynomials in $\C[x][y]$, i.e.
$\phi_k(x,y)=\phi_{k1}(x,y)\phi_{k2}(x,y)\cdots\phi_{kn_k}(x,y)$,
where each $\phi_{k\ell}(x,y)$ is an irreducible monic polynomial in
$\C[x][y]$. Therefore, every irreducible component
$\phi_{k\ell}(x,y)=0$ contains a character of the representation
corresponding to $s=i$ where $i=\sqrt{-1}$. In other words,
$\phi_{k\ell}(i+1/i,y)=0$ has a solution for $y$ for any $k,\ell$.
Since the character of this representation is 0, the representation
is irreducible metabelian (see \cite[Section 4]{NY08-1}).
\end{proof}

Let $K$ be a
2-bridge knot $K(\alpha, \beta)$ and fix a Wirtinger presentation
$G(K)=\la a,b\,|\, wa=bw\ra$. If $p=2n+1$ is a prime divisor of
$\alpha$, then as in \cite[Section 2]{HM09-1} one can define a
representation $\xi\colon G(K)\to GL(2,\C)$ by
$$
\xi(a) =
\begin{pmatrix}
-1&1\\
0&1
\end{pmatrix},\qquad
\xi(b) =
\begin{pmatrix}
-1&0\\
\omega&1
\end{pmatrix},
$$
where $\omega\in \C$ is an algebraic integer whose minimal
polynomial over $\Z$ is $z^n + c_{n-1}^{(n)}z^{n-1}+ \cdots
+c_1^{(n)}z+p$ for some $c_j^{(n)}\in \Z$, $1\le j\le n-1$.

\begin{lemma}\label{lem:4g-2}
Let $K=K(\alpha,\beta)$ be a $2$-bridge knot of genus $g$ and
$c\in\Z$ the leading coefficient of $\D_K(t)$. Let $p=2n+1$ be a
prime divisor of $\alpha$. Suppose $c\not\equiv 0\,\,
\mathrm{mod}~p$. Then the degree of $\D_{K,\xi}(t)$ is $4g-2$.
\end{lemma}

\begin{proof}
We calculate $\D_{K,\xi}(t)$ according to \cite[Section 7]{HM09-1}.
Put $r=waw^{-1}b^{-1}$.  As was shown in \cite[p.17]{HM09-1} we
obtain
$$
\Phi \left( \frac{\p r}{\p a} \right) =
\begin{pmatrix}
\D_K(-t)+\omega\mu_{11} & \mu_{12}\\
\omega\mu_{21} & \D_K(t)+\omega\mu_{11}
\end{pmatrix},
$$
where $\mu_{ij}\in \Z[\omega][t^{\pm1}]$. Then $ \D_{K,\xi}(t)=\det
\Phi \left( \frac{\p r}{\p a} \right)/(1-t^2).$ Since $K$ is a
2-bridge knot, it is alternating, and therefore $\deg(\D_K(t)) = 2g$
(see \cite{Crowell59-1,Murasugi58-1}). Therefore we can write
$\D_K(t) = ct^l + b_{l-1}t^{l-1} + \cdots + b_{m+1}t^{m+1} + ct^m$
where $l-m=2g$ and $b_j\in \Z$ $(m+1\le j \le l-1)$. Then one can
see that each $\mu_{ij}$ is written as $a_l^{ij}t^l +
a_{l-1}^{ij}t^{l-1} + \cdots + a_{m+1}^{ij}t^{m+1} + a_m^{ij}t^m$ for
some $a_k^{ij}\in\Z[\omega]$ (see \cite[Section 7]{HM09-1}).

Therefore we obtain that $\D_{K,\xi}(t) =
\psi_{2l-2}t^{2l-2}+\psi_{2l-1}t^{2l-1}+\cdots + \psi_{2m+1}t^{2m+1}
+ \psi_{2m}t^{2m}$ for some $\psi_k\in \Z[\omega]$, $2m\le k \le
2l-2$. We will show that $\psi_{2l-2}\ne 0$ and $\psi_{2m}\ne 0$,
which will complete the proof.

Since the minimal polynomial of $\omega$ over $\Z$ is $z^n +
c_{n-1}^{(n)}z^{n-1}+ \cdots +c_1^{(n)}z+p$,
$$
\psi_{2l-2}=c^2+pc_0+c_1\omega+\cdots+c_{n-1}\omega^{n-1}
$$
for some $c_j\in\Z$. If $c_j\not=0$ for some $j~(1\le j\le n-1)$,
then $\psi_{2l-2}\ne 0$. Otherwise $\psi_{2l-2}=c^2+pc$, and since
$c\not\equiv 0\,\, \mathrm{mod}~p$, $\psi_{2l-2}$ is nonzero.
Similarly, one can prove that $\psi_{2m}\ne 0$.
\end{proof}

\begin{lemma}\label{lem:nonmonic}
Let $K=K(\alpha,\beta)$ be a $2$-bridge knot and $c\in\Z$ the
leading coefficient of $\D_K(t)$. Let $p=2n+1$ be a prime divisor of
$\alpha$. Suppose
$c\not\equiv 0 \mbox{ and } c^2\not\equiv \pm1\,\, \mathrm{mod}~p$. 
Then for any irreducible metabelian
representation $\rho$, $\D_{K,\rho}(t)$ is nonmonic.
\end{lemma}

Note that $c\ne \pm 1$ by the assumption,
and therefore the knot $K$, which is alternating, is not fibered.

\begin{proof}
Fix a Wirtinger presentation $G(K)=\la a,b\,|\, wa=bw\ra$. Since
$\rho$ is irreducible metabelian, the character of $\rho$ is 0 (see
\cite[Section 4]{NY08-1}). Therefore we may assume that
$$
\rho(a) =
\begin{pmatrix}
i&1\\
0&-i
\end{pmatrix},\qquad
\rho(b) =
\begin{pmatrix}
i&0\\
2-y_0&-i
\end{pmatrix}
$$
for some $y_0\in \C$. Let $\xi\colon G(K)\to GL(2,\C)$ be the
representation defined in the paragraph preceding
Lemma~\ref{lem:4g-2}. 
Then we can easily check that $\xi$ is conjugate to $i\cdot\rho$ and
$\D_{K,\rho}(t) = \D_{K,\xi}\left(-it\right)$ holds. 
Since the highest exponent of
$\D_{K,\xi}(t)$ is $2l-2$ as we calculated in the proof of
Lemma~\ref{lem:4g-2}, this implies that if the leading coefficient
of $\D_{K,\xi}(t)$ is not $\pm1$, then $\D_{K,\rho}(t)$ is nonmonic.

By the assumption, as in the proof of Lemma~\ref{lem:4g-2} the
leading coefficient of $\D_{K,\xi}(t)$ is
$$
\psi_{2l-2}=c^2+pc_0+c_1\omega+\cdots+c_{n-1}\omega^{n-1}
$$
for some $c_j\in\Z$. If $c_j\not=0$ for some $j~(1\le j\le n-1)$,
then $\psi_{2l-2}\ne \pm1$. If $c_j=0$ for any $1\leq j\leq n-1$,
then $\psi_{2l-2}=c^2+pc_0$. Since $c^2+pc_0\equiv c^2\not\equiv
\pm1~\mathrm{mod}~p$, we have $c^2+pc_0\ne \pm1$. Therefore
$\psi_{2l-2}\ne \pm1$.
\end{proof}

Using the above lemmas, we can obtain the
following theorem which extends Theorem~\ref{thm:genus1} under a
certain condition. 

\begin{theorem}\label{thm:genus2}
Let $K=K(\alpha,\beta)$ be a $2$-bridge knot of genus $g$ and
$c\in\Z$ the leading coefficient of $\D_K(t)$. Let $p$ be an odd
prime divisor of $\alpha$.  Suppose $c\not\equiv 0 \,\,
\mathrm{mod}~p$. Then $\deg(\D_{K,\chi}(t)) = 4g-2$ for all but
finitely many $\chi\in X^\mathrm{nab}(K)$.
\end{theorem}

\begin{proof}
By Lemma~\ref{lem:metabelian}, any irreducible component of
$X^\mathrm{nab}$ contains a metabelian representation, say $\rho$.
From Lemma~\ref{lem:4g-2} and the proof of Lemma~\ref{lem:nonmonic},
one can see that $\deg(\D_{K,\rho}(t)) = 4g-2$. Following the
notations in the proof of Theorem~\ref{thm:fibered-finite-1}, for
$\chi=\chi(x,y)\in X^\mathrm{nab}$ let us write
$$
\D_{K,\chi}(t) = \psi_{4g-2}(x,y)t^{4g-2} +
\psi_{4g-3}(x,y)t^{4g-3}+\cdots +\psi_1(x,y)t +\psi_0(x,y).
$$
Then we have $\deg(\D_{K,\chi}(t)) < 4g-2$ exactly when
$\psi_{4g-2}(x,y)=0$ or $\psi_0(x,y)=0$. Since on each irreducible
component of $X^\mathrm{nab}$ there exists a (metabelian) character
whose associated twisted Alexander polynomial is of degree $4g-2$,
$\psi_{4g-2}(x,y)$ is not identically zero on 
every irreducible component
$\phi_{k\ell}(x,y)=0$. 
Similarly, $\psi_0(x,y)$ is not identically zero on $\phi_{k\ell}(x,y)=0$, 
either. Therefore by B\'{e}zout's Theorem, $\phi(x,y)=0$ and
$\psi_{4g-2}(x,y)=0$ have finitely many intersection points, and so
do $\phi(x,y)=0$ and $\psi_0(x,y)=0$.
\end{proof}

\begin{example}\label{ex:4.7}
Let $K$ be the $2$-bridge knot $K(7,3)=5_2$ with the Alexander
polynomial $\D_K(t)=2t^2-3t+2$.
This is one of twist knots and hence has genus one.
By Theorem~\ref{thm:fibered-finite-2}, the number of monic characters in
$X^\mathrm{nab}(K)$ is finite. Moreover, in this example, we can
directly find monic characters as follows. The knot group $G(K)$ has
a presentation
$$
G(K)
=
\la
a,b\,|\,wa=bw
\ra,
\quad
w=[b,a^{-1}]^2.
$$
We then see that
the polynomial $\phi(x,y)$ is given by
an irreducible polynomial
$$
\phi(x,y)
=
1-4x^2+2x^4+(2-x^2-x^4)y-(1-2x^2)y^2-y^3.
$$
The twisted Alexander polynomial in this case is
$$
\D_{K,\chi}(t)
=
\psi_2(x,y) t^2+\psi_1(x,y) t+\psi_0(x,y),
$$
where
$\psi_2(x,y)=\psi_0(x,y)=2x^2-x^2y+y^2$ and $\psi_1(x,y)=-2x$. 
An easy calculation shows that
$\{(x,y)\in {\C}^2\,|\,\phi(x,y)=0,~\psi_2(x,y)=0~\mbox{and}~\psi_1(x,y)=1\}$
is an empty set and
$$
\{(x,y)\in{\C}^2\,|\,\phi(x,y)=0~\mbox{and}~\psi_2(x,y)=1\}
=
\left(
\pm\frac{1}{\sqrt{2}},\frac12
\right),
$$
namely there are only two intersection points (see
\cite[Example 5.4]{Morifuji08-1}).

Next by Theorem~\ref{thm:genus2},
$\deg(\D_{K,\chi}(t))=4\cdot1-2=2$ for all but finitely many
$\chi\in X^{\mathrm{nab}}(K)$. In fact, the similar calculation as above
shows that
$$
\{(x,y)\in{\C}^2\,|\,\phi(x,y)=0~\mbox{and}~\psi_2(x,y)=0\}
=
\left(
\pm\frac{i}{\sqrt{6}},-\frac23
\right).
$$
Therefore except for these two points the twisted Alexander polynomial
$\D_{K,\chi}(t)$ is of degree 2.
\end{example}

\section{A bound for $2$-bridge knots}

In the previous section, we discussed finiteness of the number of
monic characters for general 2-bridge knots. However, Theorem
\ref{thm:fibered-finite-1} or Theorem \ref{thm:fibered-finite-2}
says nothing about upper bounds of the number of monic characters.
In this section, we give a bound for a family of certain 2-bridge
knots.

\subsection{A family of 2-bridge knots}
\begin{figure}[ht]
  \begin{center}
    \labellist
    \pinlabel {$k$} at 326 545
    \pinlabel {$l$} at 326 355
    \endlabellist
    \includegraphics[scale=.4]{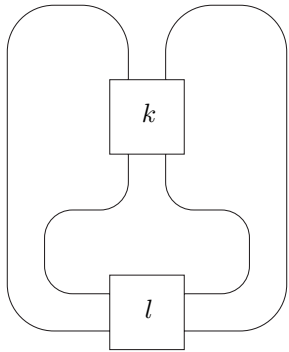}
  \end{center}
  \caption{The knot $J(k,l)$}
  \label{fig:J_kl}
\end{figure}

Let $K=J(k,l)$ be a knot as in Figure~\ref{fig:J_kl}, where $k$ and
$l$ are integers corresponding to the number of half twists in the
labeled boxes.
Positive numbers correspond to right-handed twists and negative
numbers correspond to left-handed twists respectively. The knot
$J(k,l)$ is ambient isotopic to the 2-bridge knot $K(\alpha,\beta)$
such that $\frac{\beta}{\alpha}$ is equal to $\frac{l}{1-kl}$ in
$\Q/\Z$. We remark that $J(k,l)$ is a knot if and only if $kl$ is
even (otherwise it is a 2-component link). It has an obvious
symmetry with switching $k$ and $l$. Moreover $J(-k,-l)$ is the
mirror image of $J(k,l)$. Thus, in the following, we only consider
the case where $K=J(k,2q)$ for $k>0$ and $q\in\Z$. It is easy to see
that $K=J(2,2q)$ is a twist knot and this family contains the
trefoil knot $J(2,2)$ and the figure eight knot $J(2,-2)$.

The knot group $G(K)$ for $K=J(k,2q)$
has a presentation
$$
G(K)
=
\la
a,b\,|\,(w_m)^qa=b(w_m)^q
\ra,
$$
where
$$
w_m
=
\begin{cases}
(ba^{-1})^m(b^{-1}a)^m, & k=2m;\\
(ba^{-1})^mba(b^{-1}a)^m, & k=2m+1,
\end{cases}
$$
see \cite[Proposition 1]{HS04-1}.
Let
$\rho:G(K)\to SL(2,\C)$ be a nonabelian representation defined by
the correspondence
$$
\rho(a)
=\begin{pmatrix}
s & 1 \\
0 & s^{-1}
\end{pmatrix}
\quad
\mathrm{and}
\quad
\rho(b)
=\begin{pmatrix}
s & 0 \\
-y & s^{-1}
\end{pmatrix},
$$
where
$s\not=0,y\in\C$.
For simplicity of calculations,
we have slightly modified the parametrization
of $\hat{R}^{\mathrm{nab}}(K)$.
We also denote $\rho(a),\rho(b)$ and $\rho(w_m)$
by $A,B$ and $W_m$ respectively.

The next theorem is a generalization of our previous result
(see \cite[Theorem 1.2]{GM03-1} and \cite[Theorem 5.3]{Morifuji08-1})
for the twist knot $K=J(2,2q)$.

\begin{theorem}\label{thm:5.1}
Let $K=J(k,2q)$ be a nonfibered knot,
where
$k>0,~q\in\Z$ and $k\not=2q$.
Then
the number of monic characters is bounded above by 
$$
2(k+1)^2q^2-(k+1)(k+4)|q|
$$
if $k$ is even or
$$
(k+1)(k-1)|q|
$$
if $k$ is odd.
\end{theorem}

\begin{proof}
It is known that $J(k,2q)$ is fibered only for the trivial knot
$J(k,0)$, the trefoil knot $J(2,2)$, the figure eight knot
$J(2,-2)$, the knots $J(1,2q)$ for any $q$ (it is the torus
$(2,2q-1)$-knot) and $J(3,2q)$ for $q>0$. Therefore the assertion
immediately follows from the fact that the character variety
$X^{\mathrm{nab}}(J(k,2q))$ is irreducible when $k\not=2q$ (see
\cite[Theorem 1.2]{MPV09-1}), and Propositions \ref{pro:5.2} and
\ref{pro:5.3}. That is, the upper bound is given by twice the product of
the degrees of two polynomials $\phi(x,y)$ and $\psi_{4g-2}(x,y)$.
In fact, $\{(x,y)\in{\C}^2\,|\,\phi(x,y)=0~\mbox{and}~\psi_{4g-2}(x,y)=1\}$
consists of at most
$(k+1)^2q^2-(k+1)(k+4)|q|/2$ points if $k$ is even or
$(k+1)(k-1)|q|/2$ points if $k$ is odd (B\'{e}zout's Theorem).
Similarly,
$\{(x,y)\in{\C}^2\,|\,\phi(x,y)=0~\mbox{and}~\psi_{4g-2}(x,y)=0\}$
consists of at most the same number of points,
because $\deg(\psi_{4g-2}(x,y)-1)=\deg(\psi_{4g-2}(x,y))$ clearly holds.
This completes the proof of Theorem~\ref{thm:5.1}.
\end{proof}

Since $x=s+s^{-1}$,
the Riley polynomial $\phi(s,y)\in{\Z}[s^{\pm1},y]$ and the
polynomial $\phi(x,y)$ have the same degree. Here the exponent of
the highest degree term of $\phi (s,y)$ is defined to be the degree
of $\phi(s,y)$. Moreover, we recall that the top
coefficient $\psi(s,y)$ is also a polynomial in variables
$x=s+s^{-1}$ and $y$.  Therefore $\psi(s,y)$ and $\psi(x,y)$ have the same degree.
In the following, we consider $\phi$ and
$\psi$ as polynomials in $s$ and $y$ rather than in $x$ and $y$.
The next two propositions will be proved in
Subsections 5.2 and 5.3 respectively.

\begin{proposition}\label{pro:5.2}
Let $\phi_{k,q}(s,y)$ be the Riley polynomial of $K=J(k,2q)$,
where
$k>0$ and $q\in\Z$.
Then its degree is given by
$$
\deg\phi_{k,q}(s,y)
=
\begin{cases}
(k+1)q-1, & q>0,~k\not=1; \\
(k+1)|q|, & q<0,
\end{cases}
$$
and $\deg\phi_{1,q}(s,y)=2q-2~(q>0)$.
\end{proposition}

\begin{proposition}\label{pro:5.3}
Let $\psi_{k,q}(s,y)$ be the coefficient of the highest degree term
of $\D_{K,\rho}(t)$
for $K=J(k,2q)$,
where
$k>0$ and $q\in\Z$.
Then its degree is given by
$$
\deg\psi_{k,q}(s,y)
=
\begin{cases}
(k+1)|q|- \frac{k+4}{2}, & k~is~even; \\
\frac{k-3}{2}, & k\not=1~is~odd,~q>0; \\
\frac{k-1}{2}, & k~is~odd,~q<0,
\end{cases}
$$
and $\deg\psi_{1,q}(s,y)=0$ for $q>0$.
\end{proposition}

\begin{example}\label{ex:5.3.1}
Let $K=J(4,4)=7_4$. As we have seen in Example \ref{ex:2.3},
$\phi_{4,2}(s,y)$ has the factorization $\phi_{4,2}=\phi_1\phi_2$.
Then each factor defines an irreducible curve (namely
$X^{\mathrm{nab}}(K)$ has two components, see \cite[Theorem
1.2]{MPV09-1}). So Theorem \ref{thm:5.1} does not work. Moreover
Theorem \ref{thm:fibered-finite-2} cannot be applied to $K=7_4$,
because $c=4$ and $p=3,5$ for the knot $K$, i.e. $c^2\equiv1\mod p$.
However in this example, the existence of a parabolic representation
on $\phi_2(s,y)=0$ (it is given by $\phi_2(1,y)=0$) ensures
finiteness of the number of monic characters on the whole
$X^{\mathrm{nab}}(K)$. In fact, we can show that the twisted
Alexander polynomial associated with the parabolic representation is
nonmonic.
\end{example}

\begin{remark}\label{rmk:5.3.2}
In the case of $k=2q$,
the character variety $X(J(k,k))$ is reducible and
has two components
(see \cite{MPV09-1}).
This phenomenon follows from
the fact that
the canonical component $X_0(J(k,k))$ (see Section 6) is fixed by
the involution induced by turning the 4-plat presentation upside down,
while other components are not
(see \cite{Ohtsuki94-1}).
\end{remark}

\subsection{Proof of Proposition \ref{pro:5.2}}
We first consider the case where
$K=J(k,2q)$ for $k=2m~(m>0)$ and $q>0$.

\begin{lemma}\label{lem:5.4}
For
$W_m=(BA^{-1})^m(B^{-1}A)^m=(w_m^{ij})$, we have
$$
\deg w_m^{11}=\deg w_m^{21}=2m+1~and~
\deg w_m^{12}=\deg w_m^{22}=2m.
$$
In particular,
$\deg(\tr W_m)=2m+1$.
\end{lemma}

\begin{proof}
We show them by the inductive argument on $m$.
When $m=1$, we can easily obtain
\begin{align*}
W_1
&=
BA^{-1}B^{-1}A
=
\begin{pmatrix}
1 & -s \\
-s^{-1}y & y+1
\end{pmatrix}
\begin{pmatrix}
1 & s^{-1} \\
sy & y+1
\end{pmatrix}\\
&=
\begin{pmatrix}
1-s^2y & s^{-1}-s(y+1) \\
-s^{-1}y+sy(y+1) & -s^{-2}y+(y+1)^2
\end{pmatrix}.
\end{align*}
Hence,
$\deg w_1^{11}=\deg w_1^{21}=3$ and $\deg w_1^{12}=\deg w_1^{22}=2$
as desired.
In particular,
we note that
the highest degree term in each entry $w_1^{ij}$
does not contain the negative power of $s$.

Next for the general case,
we have
\begin{align*}
W_{m+1}
&=
(BA^{-1})^{m+1}(B^{-1}A)^{m+1}
=(BA^{-1})W_m(B^{-1}A) \\
&=
\begin{pmatrix}
1 & -s \\
-s^{-1}y & y+1
\end{pmatrix}
\begin{pmatrix}
w_m^{11} & w_m^{12} \\
w_m^{21} & w_m^{22}
\end{pmatrix}
\begin{pmatrix}
1 & s^{-1} \\
sy & y+1
\end{pmatrix}
=
\begin{pmatrix}
w_{m+1}^{11} & w_{m+1}^{12} \\
w_{m+1}^{21} & w_{m+1}^{22}
\end{pmatrix},
\end{align*}
where
\begin{align*}
w_{m+1}^{11}
&=
w_m^{11}-sw_m^{21}+syw_m^{12}-s^2yw_m^{22}, \\
w_{m+1}^{12}
&=
s^{-1}w_m^{11}-w_m^{21}+(y+1)w_m^{12}-s(y+1)w_m^{22}, \\
w_{m+1}^{21}
&=
-s^{-1}yw_m^{11}+(y+1)w_m^{21}-y^2w_m^{12}+sy(y+1)w_m^{22}, \\
w_{m+1}^{22}
&=
-s^{-2}yw_m^{11}+s^{-1}(y+1)w_m^{21}-s^{-1}y(y+1)w_m^{12}+(y+1)^2w_m^{22}.
\end{align*}
Using the assumption of the induction and checking the degree of
each component carefully, we can conclude that
$\deg w_{m+1}^{11}=\deg w_{m+1}^{21}=2m+3$ and
$\deg w_{m+1}^{12}=\deg w_{m+1}^{22}=2m+2$.
In fact,
$s^2yw_m^{22}$ and
$sy(y+1)w_m^{22}$ attain the degree $2m+3$,
and
$s(y+1)w_m^{22}$ and
$(y+1)^2w_m^{22}$ attain the degree $2m+2$.
In particular,
$s^2yw_m^{22}$, which attains the degree of $\tr W_{m+1}$,
never contains the negative power of $s$.
\end{proof}

\begin{lemma}\label{lem:5.5}
For $K=J(2m,2)$,
$\deg \phi_{2m,1}(s,y)=2m$.
\end{lemma}

\begin{proof}
Using calculations in Lemma \ref{lem:5.4}, when $m=1$, it follows that
$$
\phi_{2,1}(s,y)
=
w_1^{11}+\left(s^{-1}-s\right)w_1^{12}
=
s^2+s^{-2}-y-1.
$$
Hence
we find $\deg\phi_{2,1}(s,y)=2$.
For the general case,
\begin{align*}
\phi_{2(m+1),1}
&=
w_{m+1}^{11}+\left(s^{-1}-s\right)w_{m+1}^{12}\\
&=
-sw_m^{12}+s^2w_m^{22}
+s^{-1}\left(s^{-1}w_m^{11}-w_m^{21}+(y+1)w_m^{12}-s(y+1)w_m^{22}\right),
\end{align*}
so that we have
$\deg \phi_{2(m+1),1}(s,y)=2m+2=2(m+1)$.
In fact,
$s^2w_m^{22}$ attains the degree $2m+2$.
\end{proof}

We see from \cite[Section 3]{HS04-1} that
there exists a recursive relation
$$
\phi_{k,q}
=
(\tr W_m)\phi_{k,q-1}-\phi_{k,q-2}
$$
for the Riley polynomial $\phi_{k,q}(s,y)$ of $K=J(k,2q)$.
Thus we obtain
$$
\deg\phi_{k,q}
=
\deg\phi_{k,q-1}+\deg(\tr W_m).
$$
Therefore
we see from Lemmas \ref{lem:5.4} and \ref{lem:5.5} that
the degree of $\phi_{2m,q}(s,y)$
is given by
\begin{align*}
\deg\phi_{2m,q}
&=
\deg\phi_{2m,1}+(2m+1)(q-1) \\
&=
(2m+1)q-1.
\end{align*}

Next
we consider the case where
$K=J(k,2q)$ for $k=2m+1~(m\geq 0)$ and $q>0$.
The following two lemmas are shown by the same argument
as in Lemmas \ref{lem:5.4} and \ref{lem:5.5}.

\begin{lemma}\label{lem:5.6}
For
$W_m=(BA^{-1})^mBA(B^{-1}A)^m=(w_m^{ij})$,
we have
$$
\deg w_m^{11}=\deg w_m^{21}=2m+2~and~
\deg w_m^{12}=\deg w_m^{22}=2m+1.
$$
In particular,
$\deg(\tr W_m)=2m+2$.
\end{lemma}

\begin{lemma}\label{lem:5.7}
For $K=J(2m+1,2)$,
$\deg \phi_{1,1}(s,y)=0$ and
$\deg \phi_{2m+1,1}(s,y)=2m+1~(m>0)$.
\end{lemma}

From
Lemmas \ref{lem:5.6} and \ref{lem:5.7},
we have $\deg\phi_{1,q}=2q-2$ and
\begin{align*}
\deg\phi_{2m+1,q}
&=
\deg\phi_{2m+1,1}+(2m+2)(q-1)\\
&=(2m+2)q-1.
\end{align*}

In the case of $q<0$,
we can show the following lemma.
The proof is almost the same as Lemmas \ref{lem:5.4},
\ref{lem:5.5}, \ref{lem:5.6} and \ref{lem:5.7}.

\begin{lemma}\label{lem:5.8}
$\deg\phi_{2m,-1}(s,y)=2m+1$
and
$\deg\phi_{2m+1,-1}(s,y)=2m+2$.
\end{lemma}

Since
$\deg(\tr W_m^{-1})=\deg(\tr W_m)$ holds,
we can conclude that
\begin{align*}
\deg\phi_{k,q}(s,y)
&=
\deg\phi_{k,-1}(s,y)+(k+1)(|q|-1)\\
&=
(k+1)|q|.
\end{align*}
This completes the proof of Proposition \ref{pro:5.2}.

\subsection{Proof of Proposition \ref{pro:5.3}}
Putting $r=w_m^qaw_m^{-q}b^{-1}$,
we have
\begin{align*}
\frac{\p r}{\p a}
&=
\frac{\p w_m^q}{\p a}+w_m^q
\left(
1+a\frac{\p w_m^{-q}}{\p a}
\right)\\
&=
\left(1+w_m+\cdots+w_m^{q-1}\right)\frac{\p w_m}{\p a}
+
w_m^q
\left(
1+a\left(
1+w_m^{-1}+\cdots+w_m^{-q+1}
\right)
\frac{\p w_m^{-1}}{\p a}
\right)\\
&=w_m^q\left(w_m^{-q}+\cdots+w_m^{-1}\right)\frac{\p w_m}{\p a}
+
w_m^q\left(
1-a(w_m^{-1}+\cdots+w_m^{-q})
\frac{\p w_m}{\p a}
\right)\\
&=
w_m^q\left(
1+(1-a)(w_m^{-1}+\cdots+w_m^{-q})\frac{\p w_m}{\p a}
\right),
\end{align*}
where
$$
\frac{\p w_m}{\p a}
=
-(ba^{-1}+\cdots+(ba^{-1})^m)
+
(ba^{-1})^m(1+b^{-1}a+\cdots+(b^{-1}a)^{m-1})b^{-1}
$$
if $k=2m$ or
$$
\frac{\p w_m}{\p a}
=
-(ba^{-1}+\cdots+(ba^{-1})^m)
+
(ba^{-1})^mb(1+ab^{-1}+\cdots+(ab^{-1})^m)
$$
if $k=2m+1$.

We consider the coefficient of the highest degree term of
the twisted Alexander polynomial $\D_{K,\rho}(t)$
for $K=J(2m,2q)$,
where
$m>0$ and $q>0$.
The numerator of $\D_{K,\rho}(t)$ is given by
$$
\det M_2
=
|I+(I-tA)(W_m^{-1}+W_m^{-2}+\cdots+W_m^{-q})V|,
$$
where
$$
V=-(BA^{-1}+\cdots+(BA^{-1})^m)
+(BA^{-1})^m(I+B^{-1}A+\cdots+(B^{-1}A)^{m-1})t^{-1}B^{-1}.
$$
Since
the coefficient of $t^2$ in $\det M_2$
(i.e. the highest degree term of $\D_{K,\rho}(t)$) is
$$
|-A(W_m^{-1}+\cdots+W_m^{-q})(BA^{-1}+\cdots+(BA^{-1})^m)|,
$$
we would like to compute its degree.
To this end,
we first show the following.

\begin{lemma}\label{lem:5.9}
$\deg|I+BA^{-1}+\cdots+(BA^{-1})^{m-1}|=m-1$.
\end{lemma}

\begin{proof}
When $m=1$, $\deg|I|=\deg 1=0$.
Next we assume
$$
\deg|I+BA^{-1}+\cdots+(BA^{-1})^{m-1}|
=
m-1.
$$
A direct calculation shows that
\begin{align*}
|I+BA^{-1}+\cdots+(BA^{-1})^m|
=&
|BA^{-1}+\cdots+(BA^{-1})^m|\\
&+\tr(BA^{-1}+\cdots+(BA^{-1})^m)+1\\
=&
|I+BA^{-1}+\cdots+(BA^{-1})^{m-1}|\\
&+\tr(BA^{-1}+\cdots+(BA^{-1})^m)+1,
\end{align*}
where
$BA^{-1}\in SL(2,\C)$ and
the degree of the first term is $m-1$ by the assumption
of the induction.
Further we see from the following claim that
the degree of the second term is $m$.
Therefore
we obtain
$$
\deg|I+BA^{-1}+\cdots+(BA^{-1})^m|=m.
$$
This completes the proof of Lemma \ref{lem:5.9}.
\end{proof}

\noindent
{\bf Claim.}
$\deg(\tr(BA^{-1})^m)=m$.

\begin{proof}
We set $(BA^{-1})^m=(v_m^{ij})$.
Since
$BA^{-1}
=
\begin{pmatrix}
1 & -s \\
-s^{-1}y & y+1
\end{pmatrix}$ for $m=1$,
$\deg v_1^{11}=\deg v_1^{21}=0$ and
$\deg v_1^{12}=\deg v_1^{22}=1$ hold.
Thus
$\deg(\tr BA^{-1})=1$.
In particular,
the term which attains the highest degree of
$\tr BA^{-1}$ (it is $y$ in this case)
does not contain the variable $s$.

Next
we assume that
$\deg v_m^{11}=\deg v_m^{21}=m-1$ and
$\deg v_m^{12}=\deg v_m^{22}=m$.
\begin{align*}
(BA^{-1})^{m+1}
&=
\begin{pmatrix}
v_m^{11} & v_m^{12} \\
v_m^{21} & v_m^{22}
\end{pmatrix}
\begin{pmatrix}
1 & -s \\
-s^{-1}y & y+1
\end{pmatrix}\\
&=
\begin{pmatrix}
v_m^{11}-s^{-1}yv_m^{12} & -sv_m^{11}+(y+1)v_m^{12} \\
v_m^{21}-s^{-1}yv_m^{22} & -sv_m^{21}+(y+1)v_m^{22}
\end{pmatrix}
\end{align*}
implies that
$\deg v_{m+1}^{11}=\deg v_{m+1}^{21}=m$ and
$\deg v_{m+1}^{12}=\deg v_{m+1}^{22}=m+1$.
In particular,
$yv_m^{22}$ in $(2,2)$-entry attains
the desired degree and it does not contain $s$.
Hence
$\deg(\tr (BA^{-1})^{m+1})=m+1$.
\end{proof}

\begin{lemma}\label{lem:5.10}
$\deg|I+W_m^{-1}+\cdots+W_m^{-q+1}|
=(2m+1)(q-1)$.
\end{lemma}

\begin{proof}
When $q=1$,
$\deg|I|=\deg 1=0$.
For the general case,
as in the proof of Lemma \ref{lem:5.9},
we get
\begin{align*}
|I+W_m^{-1}+\cdots+W_m^{-q+1}|
=&
|W_m^{-1}+\cdots+W_m^{-q+1}|\\
&+\tr(W_m^{-1}+\cdots+W_m^{-q+1})+1\\
=&
|I+W_m^{-1}+\cdots+W_m^{-q+2}|\\
&+\tr(W_m^{-1}+\cdots+W_m^{-q+1})+1.
\end{align*}
The degree of the first term is
$(2m+1)(q-2)$ by the assumption of the induction.
On the other hand,
we obtain
\begin{align*}
\deg\left(
\tr(W_m^{-1}+\cdots+W_m^{-q+1})
\right)
&=
\deg\left(\tr W_m^{-q+1}\right)\\
&=(q-1)\deg\left(\tr W_m^{-1}\right)\\
&=(2m+1)(q-1),
\end{align*}
because
$\deg(\tr W_m^{-1})=2m+1$ by Lemma \ref{lem:5.4},
and $\deg(\tr N^q)=q\deg(\tr N)$ holds for $N\in SL(2,\C)$ in general.
This completes the proof.
\end{proof}

Since
$A,B,W_m\in SL(2,\C)$,
Lemmas \ref{lem:5.9} and \ref{lem:5.10} show that
\begin{align*}
&\deg|-A(W_m^{-1}+\cdots+W_m^{-q})(BA^{-1}+\cdots+(BA^{-1})^m|\\
=&
\deg|I+W_m^{-1}+\cdots+W_m^{-q+1}||I+BA^{-1}+\cdots+(BA^{-1})^{m-1}|\\
=&\deg|I+W_m^{-1}+\cdots+W_m^{-q+1}|+\deg|I+BA^{-1}+\cdots+(BA^{-1})^{m-1}|\\
=&(2m+1)q-(m+2).
\end{align*}
We finished the proof for $K=J(2m,2q)$.

Next
we consider the case where $K=J(k,2q)$ for $k=2m+1~(m\geq0)$ and $q>0$.
In this case,
the numerator of $\D_{K,\rho}(t)$ is given by
$$
\det M_2
=
t^{4q}
|I+(I-tA)(t^{-2}W_m^{-1}+t^{-4}W_m^{-2}+\cdots+t^{-2q}W_m^{-q})V|,
$$
where
$$
V=-(BA^{-1}+\cdots+(BA^{-1})^m)
+t(BA^{-1})^mB(I+AB^{-1}+\cdots+(AB^{-1})^m).
$$

An easy calculation shows that
the coefficient of the highest degree term in
$\det M_2$ coincides with the top coefficient of
$$
|I-t^2A(t^{-2}W_m^{-1}+\cdots+t^{-2q}W_m^{-q})(BA^{-1})^mB
(I+AB^{-1}+\cdots+(AB^{-1})^m)|.
$$

We first consider the case where
$q=1$.

\begin{lemma}\label{lem:5.11}
The degree of
$$
|I-AW_m^{-1}
\left(BA^{-1}\right)^mB
\left(I+AB^{-1}+\cdots+(AB^{-1})^m\right)|
$$
is $m-1$ for $m>0$ and $0$ for $m=0$.
\end{lemma}

\begin{proof}
Since
$W_m^{-1}=(B^{-1}A)^{-m}A^{-1}B^{-1}(BA^{-1})^{-m}$,
we have
\begin{align*}
&|I-AW_m^{-1}
\left(BA^{-1}\right)^mB
\left(I+AB^{-1}+\cdots+(AB^{-1})^m\right)|\\
=&
|I-\left(BA^{-1}\right)^m
\left(I+AB^{-1}+\cdots+(AB^{-1})^m\right)|\\
=&
|I+BA^{-1}+\cdots+(BA^{-1})^{m-1}|.
\end{align*}
We see from Lemma \ref{lem:5.9} that
its degree is $m-1$.
When
$m=0$,
$\deg|I-AW_0^{-1}B|=\deg|O|=0$,
where
$O$ denotes the zero matrix.
\end{proof}

For the general case,
$t^{-2q}W_m^{-q}$ is newly added
under the assumption of the induction on $q$,
but
it never contributes the coefficient
of the highest degree term.
Because
it has a negative power with the variable $t$.
Namely,
the degree of the top coefficient of
$\D_{K,\rho}(t)$ is constant for the knots $K=J(2m+1,q)$
when $m$ is fixed.
Then
it is equal to $m-1=\frac{k-3}{2}~(k>1)$ or
$0~(k=1)$ by Lemma \ref{lem:5.11}.

Finally,
if we consider
$\det M_1=\Phi(\frac{\p r}{\p b})$
instead of
$\det M_2=\Phi(\frac{\p r}{\p a})$
for the case of $K=J(k,2q)$,
where $k>0$ and $q<0$,
then
we can obtain the following lemma
along the line of discussion above.

\begin{lemma}\label{lem:5.12}
The degree of the coefficient of the highest degree term
is given by
$$
\deg\psi_{k,q}(s,y)
=
\begin{cases}
(k+1)(-q)-(\frac{k}{2}+2), & k=2m,~ q<0; \\
\frac{k-1}{2}, & k=2m+1,~ q<0.
\end{cases}
$$
\end{lemma}

We omit the proof of Lemma \ref{lem:5.12} and
this completes the proof of Proposition \ref{pro:5.3}.

\section{Concluding remarks}

As is well-known, each 2-bridge knot is a torus $(2,2q-1)$-knot or a
hyperbolic knot. Hence nonfibered 2-bridge knots are all hyperbolic
knots. For a hyperbolic knot $K$ in $S^3$, it is well-known that
there is, uniquely up to conjugation, a discrete faithful
representation
$\bar{\rho}_0:G(K)\to\mathrm{Isom}^+(\mathbb{H}^3)\cong PSL(2,\C)$
such that $\mathbb{H}^3/G(K)\cong S^3-K$.
Then, due to Thurston, $\bar{\rho}_0$ can be lifted to a
discrete faithful representation $\rho_0:G(K)\to SL(2,\C)$ (see
\cite[Proposition 3.1.1]{CS83-1}). For the representation $\rho_0$,
it is natural to raise the following conjecture.

\begin{conjecture}[Dunfield-Friedl-Jackson \cite{DFJ10-1}]\label{conj:6.1} 
The twisted Alexander polynomial
$\D_{K,\rho_0}(t)$ detects all hyperbolic fibered knots.
\end{conjecture}

\begin{remark}\label{rmk:6.2}
They have shown in \cite{DFJ10-1} that Conjecture \ref{conj:6.1} is
true for all hyperbolic knots with at most thirteen crossings.
\end{remark}

Now it is easy to see that the character variety $X(K)$ always
contains a curve corresponding to abelian representations. If $K$ is
a hyperbolic knot, by Thurston's Dehn surgery theorem, there is a
so-called \textit{canonical component}\/ $X_0(K)$ in $X(K)$ which
is a curve containing the character of a discrete faithful representation $\rho_0$
(see \cite[Proposition 2]{CS84-1}).

More recently,
the work of Kronheimer-Mrowka (see also \cite{BBRW09-1})
establishes the next
general fact.

\begin{theorem}\cite{KM04-1}\label{thm:6.3}
Let $K$ be a nontrivial knot.
Then $X^{\mathrm{nab}}(K)$ contains a curve for which
all but finitely many of its elements are
characters of irreducible representations.
\end{theorem}

In our point of view to the fibering problem,
we conclude the present paper with
the following conjecture.

\begin{conjecture}\label{conj:6.4}
For a nonfibered knot $K$,
there exists a curve component $X_1(K)$ in $X^{\mathrm{nab}}(K)$
so that $\{\chi\in X_1(K)\,|\,\D_{K,\chi}(t)~\mathrm{is~monic}\}$ 
is a finite set.
\end{conjecture}

\begin{remark}\label{rmk:6.5}
If $K=J(k,l)$ is a nonfibered (hence hyperbolic) knot, the canonical
component in $X^\mathrm{nab}(K)$ has a finite number of monic
characters: if $k\ne l$, $X^{\mathrm{nab}}(K)$ is irreducible (hence
it is the canonical component) by \cite[Theorem 1.2]{MPV09-1} and it
follows from Theorem~\ref{thm:fibered-finite-1}. If $k=l$, by
\cite[Proposition 4.6]{MPV09-1} one can see that the canonical
component has a reducible representation. Thus following the
arguments in Theorems~\ref{thm:fibered-monic} and
\ref{thm:fibered-finite-1}, we see that the canonical component has
finitely many monic characters. Therefore one may further conjecture
that for a nonfibered hyperbolic knot there are only finitely many
monic characters in the canonical component. We do not know if this
conjecture is true for 2-bridge knots.
\end{remark}

\noindent \textit{Acknowledgements}. The authors would like to thank
Michel Boileau, Francis Bonahon, Teruaki Kitano, Fumikazu Nagasato
and Robert Penner for helpful comments. 
The authors also thank Stefan Friedl 
for helpful comments and correcting the proof of
Theorem~\ref{thm:fibered-finite-1}. The first author was supported
by the National Research Foundation of Korea(NRF) grant funded by
the Korea government(MEST) (No. 2009-0068877 and 2009-0086441). The
second author is supported in part by the Grant-in-Aid for
Scientific Research (No. 20740030), the Ministry of Education,
Culture, Sports, Science and Technology, Japan.


\end{document}